\documentclass[12pt,amsfonts,amscd]{amsart}

\usepackage{amssymb,amsmath,amsthm,longtable,multirow,array}
\usepackage{graphicx}

\usepackage{datetime}
\usepackage{amsfonts}
\usepackage{color}
\usepackage{mathrsfs}
\usepackage{tikz}

\newtheorem{Theorem}{Theorem}[section]
\newtheorem{Proposition}[Theorem]{Proposition}

\newtheorem{Lemma}[Theorem]{Lemma}
\newtheorem{Corollary}[Theorem]{Corollary}
\theoremstyle{definition}\newtheorem{Definition}[Theorem]{Definition}
\newtheorem{Example}[Theorem]{Example}
\newtheorem{Exercise}[Theorem]{Exercise}
\newtheorem{Question}[Theorem]{Question}
\newtheorem{Para}[Theorem]{}
\theoremstyle{remark}
\newtheorem*{remark*}{Remark}

\providecommand\ba[1]{\begin{align*}#1\end{align*}}

\providecommand\baaa[1]{\begin{equation}\begin{split}#1\end{split}\end{equation}}

\providecommand\brs{\begin{remark*}}
\providecommand\ers{\end{remark*}}

\providecommand\be{\begin{enumerate}}
\providecommand\ee{\end{enumerate}}
\providecommand\bT{\begin{Theorem}}
\providecommand\eT{\end{Theorem}}
\providecommand\bP{\begin{Proposition}}
\providecommand\eP{\end{Proposition}}
\providecommand\bD{\begin{Definition}}
\providecommand\eD{\end{Definition}}
\providecommand\bE{\begin{Example}}
\providecommand\eE{\end{Example}}
\providecommand\bEE{\begin{Exercise}}
\providecommand\eEE{\end{Exercise}}
\providecommand\bPP{\begin{Para}}
\providecommand\ePP{\end{Para}}
\providecommand\bL{\begin{Lemma}}
\providecommand\eL{\end{Lemma}}
\providecommand\bC{\begin{Corollary}}
\providecommand\eC{\end{Corollary}}
\providecommand\bpp{\begin{proof}}
\providecommand\epp{\end{proof}}
\providecommand\bee{\begin{equation}}
\providecommand\eee{\end{equation}}
\providecommand\bQ{\begin{Question}}
\providecommand\eQ{\end{Question}}

\providecommand\beqq{\begin{eqnarray*}}
\providecommand\eeqq{\end{eqnarray*}}

\providecommand\bay{\begin{array}}
\providecommand\eay{\end{array}}

\providecommand\CC{{\Bbb C}}
\providecommand\RR{{\Bbb R}}

\providecommand\NN{{\Bbb N}}

\providecommand\DD{{\Bbb D}}

\providecommand\TT{{\Bbb T}}

\providecommand\pd{\partial}

\providecommand\al{\alpha}
\providecommand\bt{\beta}

\providecommand\dl{\delta}

\providecommand\eps{\varepsilon}
\providecommand\Gm{\Gamma}
\providecommand\gm{\gamma}

\providecommand\vp{\varphi}
\providecommand\Lm{\Lambda}
\providecommand\lm{\lambda}

\providecommand\Om{\Omega}
\providecommand\om{\omega}

\providecommand\zt{\zeta}

\providecommand\sbs{\subset}

\providecommand\arg{\operatorname{arg}}

\providecommand\diam{\operatorname{diam}}

\providecommand\rep {\operatorname{Re}}
\providecommand\imp {\operatorname{Im}}

\providecommand\qs{quasisymmetric\ }
\providecommand\qst{quasisymmetric}

\providecommand\qc{quasicircle\ }
\providecommand\qct{quasicircle}
\providecommand\qcs{quasicircles\ }
\providecommand\qcst{quasicircles}

\providecommand\qd{quasidisk\ }
\providecommand\qdt{quasidisk}
\providecommand\qds{quasidisks\ }
\providecommand\qdst{quasidisks}

\providecommand\iiff{if and only if }

\begin{document}

\title[$t$-quasi-circles]{Analytic characterization of $t$-quasicircles and conformal mappings onto $t$-quasidisks}

\author{Xin Wei}

\begin{abstract} 
In this paper, we introduce a new class of mappings defined on the unit circle $\TT$, termed $(\rho,t)$-\qs mappings, which generalizes the classical concept of \qs mappings. Using this broader class of mappings, we provide an analytic characterization of $t$-\qcst. This result can be viewed as a $t$-\qc analogue of a fundamental theorem by Tukia and V\"ais\"al\"a \cite{TV}. Furthermore, we study conformal mappings from the unit disk $\DD$ onto $t$-\qds and show that their boundary values are $(\rho,t^6)$-\qst. This result generalizes the Quasicircle Theorem to the case of $t$-\qcst.
\end{abstract}


\address{xwei@xsyu.edu.cn, School of Science, Xi'an Shiyou University, Xi'an, Shaanxi 710065, P.R.China}

\maketitle

\bigskip
\section{Introduction}\label{Intro}
A Jordan curve $\Gm$ is a \qc if it is the image of the unit circle under a quasiconformal mapping of the extended complex plane onto itself. The bounded component of $\RR^2\setminus\Gm$ is called a \qdt. Quasicircles appear in many different settings in complex analysis and geometric function theory. For example, in complex dynamics, they are the Julia sets of some rational mappings; also \qcs play an important role in universal Teichm\"uller space, through \qs homeomorphism of the unit circle. Besides the original definition, Ahlfors  established a geometric characterization of \qcs in \cite{A}, which is also called bounded turning condition. Let $\Gm$ be a Jordan curve, for any $ x,y\in \Gm$, let $\Gm(x,y)$ be the subarc of $\Gm$ between $x$ and $y$ which has a smaller diameter, or, be either subarc when both have the same diameter. The diameter of $\Gm(x,y)$ is denoted by $\diam\Gm(x,y)$.

\bT(Ahlfors' Condition)
Jordan curve $\Gm$ is a \qc \iiff there exists a constant $C\ge 1$ such that 
\baaa{\label{2509031}\diam\Gm(x,y)\le C|x-y|}
for any $x,y\in \Gm$.
\eT

\begin{figure}[htbp]
\includegraphics[width=0.2\textwidth]{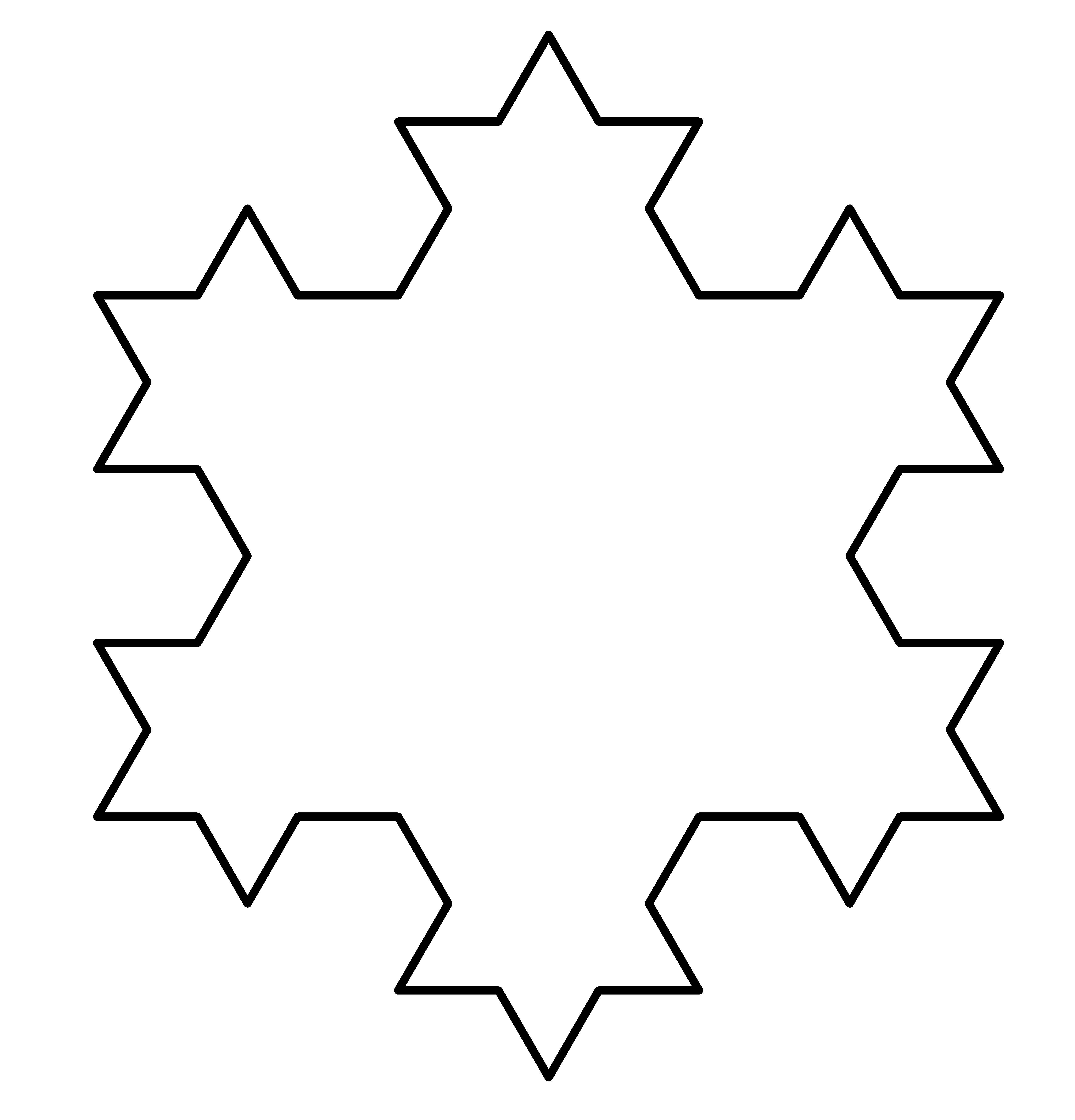}
\caption{The third step of the construction of Koch snowflake curve}\label{koch}
\end{figure}

\bE
 The Koch snowflake curve is a \qct; see Figure~\ref{koch}. 
\eE

This equivalent characterization of \qct s is also applied to define quasiarcs, i.e., a planar arc $\Gm$ is a quasiarc \iiff there exists a constant $C\ge 1$ such that (\ref{2509031}) is satisfied for any $x,y\in \Gm$. Here $\Gm(x,y)$ is the unique subarc of $\Gm$ connecting $x$ and $y$. In \cite{N}, in order to study the Whitney phenomenon, Norton introduced the definition of $t$-quasiarcs for $t\ge 1$ .

\textit{ A planar arc $\Gm$ is a $t$-quasiarc \iiff there exists a constant $C\ge 1$ such that $$(\diam\Gm(x,y))^t\le C|x-y|$$ for any $x,y\in \Gm$. }
 
We extend Norton's construction to establish the definition of $t$-\qcs for $0<t\le 1$. The interior bounded domain of a $t$-\qc is called a $t$-\qdt. 

\bD\label{2503311}
Let $0<t\le 1$, and let $\Gm\sbs\RR^2$ be a Jordan curve. Then $\Gm$ is a $t$-\qc \iiff there exists a constant $C\ge 1$ such that 
\ba{\diam\Gm(x,y)\le C|x-y|^t}
for any $x,y\in \Gm$.  In particular, a $1$-\qc is a usual \qct.
\eD

\bE
In Figure~\ref{onehalf}, the arc $AB$ is a circular arc and the arc $CB$ is a straight line tangents $AB$ at $B$. According to the definition of $t$-\qcst, the arc $ABC$ is a $1/2$-quasiarc.
\eE

In \cite{MWW}, the authors studied the local properties of self-similar $t$-\qcst. Also, they provided a method for constructing $t$-\qcs for any given $t\in(0,1]$. Also, a global flatness condition of $t$-\qcs was studied in \cite{WW}, and it follows that the constant distance boundary of a $t$-\qc is also a $t$-\qct.

In \cite{TV}, Tukia and V\"ais\"al\"a provided another analytic  characterization of \qcs of metric spaces by using \qs mappings. We provide the definition of \qs mappings between metric spaces firstly. In the following definition, the notation $|a-b|$ stands for the distance between $a$ and $b$ in any metric space.

\bD\label{2503131}\cite{TV}
Let $f: X\to Y$ be a homeomorphism between metric spaces $X$ and $Y$, and let $\rho: [0,\infty)\to [0,\infty)$ be a strictly increasing continuous function. The mapping $f$ is called \qs or $\rho$-\qs if $|a-x|\le k|b-x|$ implies $$|f(a)-f(x)|\le \rho(k) |f(b)-f(x)|$$ for any $a,b,x\in X$ and $k\ge 0$. While, a homeomorphism $f: X\to Y$ is called weakly \qs or weakly $R$-\qs if for any $a,b,x\in X$, $|a-x|\le |b-x|$ implies $$|f(a)-f(x)|\le R |f(b)-f(x)|$$ for some constant $R$. 
\eD

According to the next theorem, we see that \qcs and \qs mappings are closely related.  

\bT\cite[Theorem~4.9]{TV}
A curve is a \qc \iiff 
it is the image of the unit circle under a \qs mapping.
\eT 

We would like to create a parallel relationship for $t$-\qcst. So, in Section~\ref{anf}, we generalize the definition of weakly $R$-\qs mappings and obtain weakly $(R,t)$-\qs mappings (see Definition~\ref{23053001}). Moreover, we prove that
\textit{any $t$-\qc is the image of the unit circle under a weakly $(R,t)$-\qs mapping} (see Theorem~\ref{2506062}). 
This result provides us an analytic characterization of $t$-\qcst. We also introduce the definition of $(\rho,t)$-\qs mappings in this section, and establish the relations between weakly $(R,t)$-\qs and $(\rho,t)$-\qs mappings. 

Let $\Om$ be a \qdt, which is the inner domain of a \qct, and then it is a bounded simply connected domain. The Riemann Mapping Theorem guarantees that there exists a conformal mapping $f$ such that $f(\DD)=\Om$. The boundary behaviour of $f$ was studied in \cite{P}, which is called  Quasicircle Theorem (see page 94 of \cite{P}). 

\bT\cite[Quasicircle Theorem]{P} Let $\Gm$ be a Jordan curve in $\CC$ and let $f$ map $\DD$ conformally onto the inner domain $\Om$ of $\Gm$. Then the following conditions are equivalent.

(1) $\Gm$ is a \qct;

(2) $f$ is \qs on $\TT$.
\eT

In Section~\ref{pocm}, we introduce several auxiliary results on conformal mappings. Moreover, some estimates about $|f'|$ are obtained, which are fundamental to the study of boundary behaviour of conformal mappings which map the unit disk onto $t$-\qdst. 

In Section~\ref{cmotd}, the boundary behaviour of conformal mappings from $\DD$ onto $t$-\qds are studied. We obtain, in Theorem~\ref{2512034}, that
\textit{the boundary value of conformal mapping $f$ from $\DD$ onto a $t$-\qd $\Om$ is $(\rho,t^6)$-\qst.} In addition, we can replicate the Quasicircle Theorem from our result by letting $t=1$.

\brs
The open ball $\{y\in \CC: |x-y|<r\}$ is written as $\DD(x,r)$ and the closed ball $\{y\in \CC: |y-x|\le r\}$ is written as $\bar \DD(x,r)$. Specially, the open and closed unit disks are denoted by $\DD$ and $\bar\DD$ respectively. The unit circle is denoted by $\TT$. For any Jordan curve $\Gm$, the notation $\Gm(x,y)$ represents the subarc of $\Gm$ between $x$ and $y$ with a smaller diameter, or, represents either subarc when both have the same diameter; and $\diam\Gm(x,y)$ represents its diameter.
\ers

\section{Generalized \qs mappings}\label{anf}

In \cite{TV}, it was proved that any \qc is the image of the unit circle $\TT$ under a \qs mapping. So it is natural to find this kind of analytic characterization of $t$-\qcst. Thus, in this section, we define a class of mappings between metric spaces which can be viewed as a generalization of the  \qs mappings of Definition \ref{2503131}. And then the relationships between this kind of mappings and $t$-\qcs are studied; the main result states that \textit{every $t$-\qc is the image of the unit circle under a mapping we defined in this section.}

Let $0<t\le 1$, and let $X$, $Y$ be metric spaces. We use the notation $|a-b|$ to denote the distance between $a$ and $b$ in any metric space.  

\bD\label{23053001}
Let $f: X \to Y$ be a homeomorphism, and let $\rho: [0,\infty)\to [0,\infty)$ be a strictly increasing continuous function. The mapping $f$ is called $(\rho,t)$-\qs if $|a-x|\le k|b-x|$ implies 
\ba{|f(a)-f(x)|\le \rho(k) |f(b)-f(x)|^t} 
for any $a,b,x\in X$ and $k\ge 0$. While, say $f$ is weakly $(R,t)$-\qs if $|a-x|\le |b-x|$ implies 
\ba{|f(a)-f(x)|\le R|f(b)-f(x)|^t}
for some constant $R>0$.
\eD

When $t=1$, a $(\rho,t)$-\qs mapping is reduced to a $\rho$-\qs mapping and a weakly $(R,t)$-\qs mapping is reduced to a weakly $R$-\qs mapping. Also, directly follows  from the definitions, every $(\rho,t)$-\qs mapping is weakly $(R,t)$-\qs with $R=\rho(1)$.  When $t=1$, every weakly $R$-\qs mapping is $\rho$-\qs (see Theorem~ 2.16 of \cite{TV}). However, in the case where $0<t<1$, whether an analogous conclusion holds remains an open question.

\bE
If $f: X \to Y$ satisfies the bi-H\"older condition, i.e., there exist positive constants $K_1,K_2$ such that
\ba{K_1 |x-y|^{1/t}\le |f(x)-f(y)|\le K_2 |x-y|^t}
for any $x,y\in X$, then $f$ is weakly $(R,t^2)$-\qst.
\eE
\bpp
For any $a,b,x\in X$ with $|a-x|\le |b-x|$, then $|b-x|^{1/t}\le 1/K_1|f(b)-f(x)|$. Thus $|a-x|^t\le 1/K_1^{t^2}|f(b)-f(x)|^{t^2}$. Since $|f(a)-f(x)|\le K_2|a-x|^t$, it leads to $|f(a)-f(x)|\le K_2/K_1^{t^2}|f(b)-f(x)|^{t^2}$. Let $R=K_2/K_1^{t^2}$. Then the proof is complete.
\epp

The next proposition shows that the image of the unit circle $\TT$ under a weakly $(R,t)$-\qs mapping is a $t$-\qct. This is our first conclusion about the connection between weakly $(R,t)$-\qs mappings and $t$-\qcst. 

\bP\label{2506051}
Let $f:\TT\to \RR^2$ be a weakly $(R,t)$-\qs mapping. Then $\Gm:=f(\TT)$ is a $t$-\qct.
\eP
\bpp
For any $a,b\in \Gm$, there exist $x,y\in \TT$ such that $f(x)=a$ and $f(y)=b$. Let $\TT(x,y)$ be the subarc of $\TT$ between $x$ and $y$ with a smaller diameter. Then we obtain that $\diam \TT(x,y)=|x-y|$. For any $z\in \TT(x,y)$, it follows that $|z-x|\le \diam \TT(x,y)= |x-y|$. The weakly $(R,t)$-\qs property of $f$ yields that
\ba{|f(z)-f(x)|\le R|f(x)-f(y)|^t.} 
Since there are $f(u), f(v)\in f(\TT(x,y))$ such that $\diam f(\TT(x,y))=|f(u)-f(v)|$, we know that 
\ba{\diam f(\TT(x,y))=|f(u)-f(v)|&\le |f(u)-f(x)|+|f(v)-f(x)|\\
&\le 2R|f(x)-f(y)|^t.} 
Because $f(\TT(x,y))$ is a subarc of $\Gm$ between $a$ and $b$, and $\Gm(a,b)$ stands for the subarc of $\Gm$ between $a$ and $b$ which has a smaller diameter. It follows that 
\ba{\diam \Gm(a,b)\le \diam f(\TT(x,y))\le 2R|a-b|^t.}
Therefore, by its definition, $\Gm$ is a $t$-\qct.
\epp

An analogous proof yields the following proposition, which is a special case of Theorem~2.11 of \cite{TV} when $t=1$. 

\bP\label{2503141}
Let $\Lm$ be a \qct. If $f:\Lm\to \RR^2$ is $(\rho,t)$-\qs then $L:=f(\Lm)$ is a $t$-\qct. Particularly, $L$ is a \qc if $t=1$.
\eP

The above two propositions established that a weakly $(R,t)$-\qs mapping maps $\TT$ onto a $t$-\qct, and a $(\rho,t)$-\qs mapping maps \qcs onto $t$-\qcst. Next we derive the converse relation, i.e.,  any $t$-\qc is the image of $\TT$ under a weakly $(R,t)$-\qs mapping. 

Our proof is parallel to that of Theorem~4.9 of \cite{TV}, so we should introduce a definition of \cite{TV} which will be used later. A metric space $X$ is called homogeneously totally bounded, if there is an increasing function $h:[1/2,\infty)\to [1,\infty)$ such that for every $\al\ge 1/2$, every closed ball $\{y\in X: |y-x|\le r\}$ of $X$ can be covered with sets $A_1,\dots, A_s$ such that $s\le h(\al)$ and $\diam(A_j)<r/\al$ for $j=1,\dots, s$. We care about the  special case when $X\sbs\RR^2$, it follows from Remarks~2.8 of \cite{TV} that every subset of $\RR^2$ is homogeneously totally bounded with $h(\al)=4(\al \sqrt{2}+1)^2$.

The next lemma generalizes Lemma~4.4 of \cite{TV} to the case of $t$-\qcst. In the following lemma, we assume that the $t$-\qc $\Gm$ is oriented, and any  subarc has the induced orientation. For a subarc $L\sbs \Gm$, denoted by $A(L)$ the initial point of $L$, meanwhile $B(L)$ denotes the terminal point of $L$.

\bL\label{2505151}
Let $\Gm$ be a $t$-\qct. Then for any $0<\eps<1$ there exist $\dl>0$ and $p\in \NN$ such that there is a subdivision $\Gm_1,\dots \Gm_p$ of $\Gm$ with 
\ba{\dl\diam(\Gm)\le \diam(\Gm_i)\le \eps \diam(\Gm)}
for $i=1,\dots, p$.
\eL
\bpp
By a change of scale we assume that $\diam(\Gm)=1$. Select an arbitrary point $a\in \Gm$. Since $\eps<1$, there exists a subarc $\Gm_1\sbs \Gm$ which is the maximal subarc of $\Gm$ such that $A(\Gm_1)=a$ and that $\Gm_1\sbs \bar \DD(a,\eps/4)$. So $\eps/4\le \diam(\Gm_1)\le \eps/2$. Let $\Gm_2$ be the maximal subarc of $\Gm\setminus \Gm_1$ such that $A(\Gm_2)=B(\Gm_1)$ and that $\Gm_2\sbs \bar \DD(B(\Gm_1),\eps/4)$. Suppose that $\Gm_1,\dots,\Gm_j$ have been inductively constructed. If $\Gm_1\cup \cdots\cup \Gm_j=\Gm$, let $p:=j$ and the process ends. Otherwise let $\Gm_{j+1}$ be the maximal subarc of $\Gm\setminus (\Gm_1\cup \cdots\cup \Gm_j)$ such that $A(\Gm_{j+1})=B(\Gm_j)$ and that $\Gm_{j+1}\sbs \bar \DD(B(\Gm_j),\eps/4)$. Because of the compactness of $\Gm$, we could obtain subarcs $\Gm_1,\dots \Gm_{n+1}$ of $\Gm$ such that $\Gm_1\cup \cdots\cup \Gm_{n+1}=\Gm$ and that $\diam(\Gm_i)\le \eps/2$ for all $i\le n+1$ and $\diam(\Gm_i)\ge \eps/4$ for $i\le n$. Deleting $\Gm_{n+1}$ and replacing by $\Gm_n$ by the subarc with the initial point $A(\Gm_n)$ and the terminal point $a$, we get a subdivision $\Gm_1,\dots,\Gm_n$ of $\Gm$ such that $\eps/4\le \diam(\Gm_i)\le \eps$ for all $i\le n$.

Assume that $1\le i<j\le n$ and let $\Gm'$ be a subarc of $\Gm$ with endpoints $A(\Gm_i)$ and $A(\Gm_j)$. Then $\Gm_i\sbs\Gm'$ and then $\diam(\Gm')\ge \eps/4$. While $\Gm''=\Gm\setminus\Gm'$ is also a subarc between $A(\Gm_i)$ and $A(\Gm_j)$. It is not hard to see that $\diam(\Gm'')\ge \eps/4$, since $\Gm_n\sbs\Gm''$. The $t$-\qc property of $\Gm$ yields that
\ba{\eps/4\le C|A(\Gm_i)-A(\Gm_j)|^t.}
Here $C\ge 1$ is a suitable constant. It follows that $|A(\Gm_i)-A(\Gm_j)|\ge (\eps/4C)^{1/t}$. By Remark 2.8 of \cite{TV}, $\Gm$ is $h$-homogeneously totally bounded with $h(\al)=4(\al \sqrt{2}+1)^2$. Thus $n\le h((4C/\eps)^{1/t})=:p$, where $4C/\eps> 1$ since $0<t, \eps<1$ and $C\ge 1$. If $n=p$, there is nothing more to be done. Otherwise we divide $\Gm_1$ into two subarcs of diameter at least $\eps/8$. Repeating this $p-n$ times we obtain a new subdivision written as $\Gm_1,\dots, \Gm_p$ such that 
\ba{2^{-p-1}\eps\le \diam(\Gm_i)\le \eps}
for all $i\le p$. Let $p=h((4C/\eps)^{1/t})$ and $\dl=2^{-p-1}\eps$, and then the proof is complete.
\epp

With Lemma~\ref{2505151}, the next lemma follows by adapting the proof of Lemma~4.7 of \cite{TV}; we omit the details. However, some notations should be introduced. As in Lemma~\ref{2505151}, let $p=h((8C)^{1/t})$, $\dl=2^{-p-2}$. Define $\mu=2^{p+1}$, then $\mu\ge 2$ and $\mu^2\ge 2^{p+2}$. Choose positive integers $m\ge 2$ and $n$ such that 
\ba{2^n\ge \mu^m,\\
4n+3\le p^m.}
To see the existence of such $m$ and $n$, find $m$ and $n$ that satisfy the first inequality firstly. If they do not satisfy the second inequality, double them, and repeat this until the second inequality is also satisfied. Let $N=p^m+2n(p-1)$. Define $W$ to be the set of words which consist of finite sequence $w=n_1\cdots n_q$ where $1\le n_i \le N$ and $q\ge 0$; also the empty sequence $\emptyset$ is in $W$. The length of a word $w=n_1\cdots n_q$, denoted by $\ell(w)$, is the number $q$. Two words $w=n_1\cdots n_q$ and $v=m_1\cdots m_r$ can be concatenated to a new word $wv=n_1\cdots n_q m_1\cdots m_r$.

\bL\label{23061401}
Let $\Gm$ be a $t$-\qct. Then there are subarcs $\Gm_w\sbs \Gm$, $w\in W$ such that:

 (1) $\Gm_\emptyset=\Gm$.
 
 (2) $\Gm_{w 1},\dots,\Gm_{w N}$ is a subdivision of $\Gm_w$ for $w\in W$.
 
 (3) $\diam(\Gm_w)\le 2^{-\ell(w)}\diam(\Gm)$.
 
 (4) If $w, v\in W$, $\ell(w)=\ell(v)$, and if $\Gm_w$ and $\Gm_{v}$ are adjacent, then 
 \ba{\mu^{-m-2}\le \diam(\Gm_w)/\diam(\Gm_v)\le \mu^{m+2}.}
\eL

We now deduce the principal theorem of this section from Lemmas~\ref{2505151}-\ref{23061401}:

\bT\label{2506062}
If $\Gm$ is a $t$-\qc then there exists a weakly $(R,t)$-\qs mapping $f:\TT\to\RR^2$ such that $\Gm=f(\TT)$.
\eT
\bpp
For any $a,b\in \TT$, recall that $\TT(a,b)\sbs\TT$ is the subarc between $a$ and $b$ with a smaller diameter. Let $|a-b|$ denote the $1/2\pi$ of the arclength of $\TT(a,b)$. Then this is a comparable metric with the euclidean distance of $a$ and $b$ in $\RR^2$. 

Let $N$ be the constant which is defined before Lemma~\ref{23061401}. Divide $\Gm$ into $N$ subarcs $\{\Gm_w\}$ with specific properties of Lemma~\ref{23061401}, and equally divide $\TT$ into $N$ subarcs $\{I_w\}$. There is a homeomorphism $f$ such that $f(I_w)=\Gm_w$ for all $w=1,\dots, N$. We show that $f$ is weakly $(R,t)$-\qst. Without loss of generality, we assume $\diam(\Gm)=1$ by scaling. For any distinct $a,b,x\in \TT$ with $|a-x|\le |b-x|$, let $I:=\TT(b,x)$, while $J:=\TT\setminus I$.

Case 1. Assume that $f(I)=\Gm(f(b),f(x))$. Let $s$ be the smallest integer such that $N^{-s}\le |b-x|/2$. Then $I$ contains an arc $I_\al$ with $\ell(\al)=s$. It follows that
\baaa{\label{23061402}\diam(\Gm_\al)=\diam(f(I_\al)) \le \diam(f(I))\le C |f(b)-f(x)|^t.}
Here, $C$ is the constant in the definition of $t$-\qct. Meanwhile, since $|a-x|/2\le |b-x|/2\le N^{-s+1}$ it follows that $|a-x|\le 2N^{-s+1}$ and $|a-b|\le |a-x|+|b-x|\le 4N^{-s+1}$. Hence there is a sequence of arcs $I_{\al(1)},\dots,I_{\al(j)}$ following one another such that $\ell(\al(i))=s$, $j\le 4N+1$, $\al(1)=\al$ and $a\in I_{\al(j)}$. In the same manner, $|b-x|\le 2N^{-s+1}$ implies that there is a sequence of arcs $I_{\bt(1)},\dots,I_{\bt(j')}$ such that $\bt(i)\in W$, $\ell(\bt(i))=s$, $j'\le 2N+1$, $\bt(1)=\al$ and $x\in \bt(j')$. Let $c:=\mu^{m+2}$, the notation from Lemma~\ref{23061401}. Then it implies that for every $u\in I_\al$ we obtain that 
\ba{|f(a)-f(u)|\le (1+c+\cdots+c^j)\diam(\Gm_\al)\le 5Nc^{5N}\diam(\Gm_\al)} and 
\ba{|f(x)-f(u)|\le (1+c+\cdots+c^{j'})\diam(\Gm_\al)\le 3Nc^{3N}\diam(\Gm_\al).}
Combining these two inequalities with (\ref{23061402}) establish  
\baaa{\label{2505222} |f(a)-f(x)|\le 8Nc^{5N}\diam(\Gm_\al)\le 8Nc^{5N} C|f(b)-f(x)|^t.}

Case 2. If $f(J)=\Gm(f(b),f(x))$, assume that $s$ is the smallest integer such that $N^{-s}\le |b-x|/2$. Then $J$ contains an arc $J_\al$ with $\ell(\al)=s$. Then 
\baaa{\label{2505221}\diam(\Gm_\al)=\diam(f(J_\al)) \le \diam(f(J))\le C |f(b)-f(x)|^t.}
Since $|b-x|\le 2N^{-s+1}$, there exist at most $2N+2$ arcs $I_w$ with $\ell(w)=s$ such that $I_w\cap I\not=\emptyset$. Also, $|a-x|\le 2N^{-s+1}$ implies that $|a-b|\le |a-x|+|b-x|\le 4N^{-s+1}$. If we select $J_\al$ close to $b$ enough, there exist a sequence of arcs $J_{\al(1)},\dots,J_{\al(j)}$ following one another such that $\ell(\al(i))=s$, $j\le 4N+2$, $\al(1)=\al$ and $a\in J_{\al(j)}$. While there is a sequence of arcs $J_{\bt(1)},\dots,J_{\bt(j')}$ such that $\bt(i)\in W$, $\ell(\bt(i))=s$, $j'\le 2N+2$, $\bt(1)=\al$ and $x\in \bt(j')$. Then we claim that (\ref{2505221}) also implies (\ref{2505222}) by a similar argument to that of Case 1.

Combining the conclusions of Case 1 and Case 2, we prove that $f$ is weakly $(R,t)$-\qs with $R=8Nc^{5N} C$.
\epp

\section{Preparations of conformal mappings on $\DD$}\label{pocm}

In this section some auxiliary results on conformal mappings are introduced. Basically, they concern about the estimate of $|f'(z)|$ for conformal mapping $f:\DD\to \CC$. For the most part, these results are well-known; however, we  establish a more accurate estimate for the case of $|z|\ge 1/2$. Through this section, let $M_1,\;M_2\dots$ denote suitable constants. 

Recall that if $f=u+iv$ is holomorphic on $\DD$, it satisfies the Cauchy-Riemann equations on $\DD$. If we rewrite these equations in polar coordinates, it follows that 
\ba{\frac{\pd u}{\pd r}=\frac{1}{r}\frac{\pd v}{\pd \theta},\;\; \frac{\pd v}{\pd r}=-\frac{1}{r}\frac{\pd u}{\pd \theta}.}
Suppose that $f:\DD\to\CC$ is conformal, thus $\log f'=\log|f'|+i\arg f'$ is holomorphic on $\DD$. So we know that
\baaa{\label{2505273}\frac{\pd }{\pd \theta}\log |f'(z)|=-r\frac{\pd}{\pd r} \arg f'(z).}
Let $z=re^{i\theta}$. Then 
\ba{z\frac{f''(z)}{f'(z)}=z\frac{\pd}{\pd z}\log f'(z)=r\frac{\pd}{\pd r}\log f'(z)=r\frac{\pd }{\pd r}\log|f'(z)|+ir\frac{\pd }{\pd r}\arg f'(z).}
By (\ref{2505273}), we find the imaginary part
\ba{ \imp\frac{zf''(z)}{f'(z)}=r\frac{\pd}{\pd r}\arg f'(z)=-\frac{\pd}{\pd\theta}\log |f'(z)|.}
The Koebe Distortion Theorem states that
\ba{|\frac{zf''(z)}{f'(z)}-\frac{2r^2}{1-r^2}|\le \frac{4r}{1-r^2}.}
It implies that
\ba{\frac{2r^2-4r}{1-r^2}\le \rep \frac{zf''(z)}{f'(z)}\le \frac{2r^2+4r}{1-r^2},}
which is 
\baaa{\label{2506041}\frac{2r-4}{1-r^2}\le \frac{\pd}{\pd r}\log |f'(z)| \le \frac{2r+4}{1-r^2}.}
Meanwhile, we also have
\ba{-\frac{4r}{1-r^2}\le \imp \frac{zf''(z)}{f'(z)}\le \frac{4r}{1-r^2},}
which is 
\baaa{\label{2505231}-\frac{4r}{1-r^2}\le \frac{\pd}{\pd \theta}\log |f'(z)| \le \frac{4r}{1-r^2}.}

Let $\zt_1=r_1e^{i\al}$, $\zt_2=r_2e^{i\bt}$ where $0<r_1\le r_2<1$ and $\al,\bt\in (-\pi,\pi]$. Define $\zt=r_2e^{i\al}$. Let $\gm_1:=\{r_2e^{i\theta}|\; \theta: \bt\to\al \}$ and $\gm_2:=\{re^{i\al}|\; r: r_2\to r_1\}$. Here $\gm_1$ is the circular curve from $\zt_2$ to $\zt$ with radius $r_2$ and $\gm_2$ is the line segment from $\zt$ to $\zt_1$ with argument $\al$. Then $\gm:=\gm_1\cup\gm_2$ is a piecewise smooth curve from $\zt_2$ to $\zt_1$. Then (\ref{2506041}) and (\ref{2505231}) imply that
\baaa{\label{2505241}\log\frac{|f'(\zt_1)|}{|f'(\zt_2)|}=&\int_\gm d\log |f'(z)|=\int_\gm  \frac{\pd}{\pd r}\log |f'(z)| dr +\frac{\pd}{\pd \theta}\log |f'(z)| d\theta\\
=&\int_{\gm_1} \frac{\pd}{\pd \theta}\log |f'(z)| d\theta + \int_{\gm_2} \frac{\pd}{\pd r}\log |f'(z)| dr\\
\le & \frac{4r_2}{1-r_2^2}|\al-\bt|+ \int_{r_2}^{r_1} \frac{2r-4}{1-r^2} dr\\
=& \frac{4r_2}{1-r_2^2}|\al-\bt|+\log\frac{1-r_1}{1-r_2}+3\log\frac{1+r_2}{1+r_1}\\
\le & \frac{4r_2}{1-r_2^2}|\al-\bt|+\log\frac{1-r_1}{1-r_2}+\log 8 .}
By (\ref{2505241}), an estimate about $|f'(\zt_1)|/|f'(\zt_2)|$ is established, which is crucial in the proof of Theorem~\ref{2505274}.

\begin{figure}[htbp]
\includegraphics[width=0.6\textwidth]{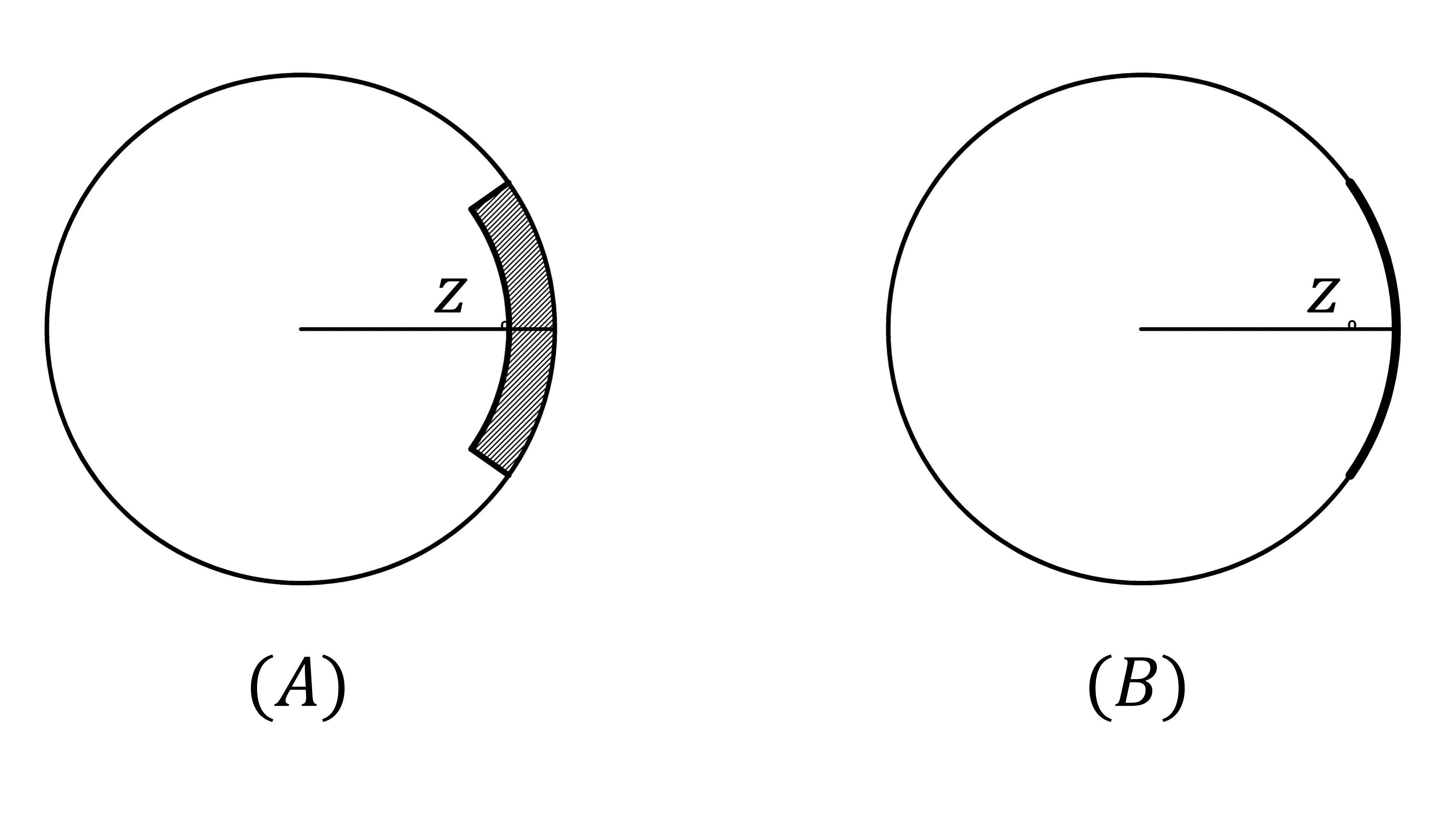}
\caption{The shaded domain in (A) is $B(z)$; the bold arc in (B) is $I(z)$}\label{BIz}
\end{figure}

For any complex number $re^{i\al}\in\DD$, we define
\ba{B(re^{i\al}):=\{\rho e^{i\theta}: r\le \rho\le 1,\; |\theta-\al|\le \pi(1-r) \},}
and 
\ba{ I(re^{i\al}):=B(re^{it})\cap \TT=\{e^{i\theta}: |\theta-\al|\le \pi(1-r)\}.}
See Figure~\ref{BIz} for examples of $B(z)$ and $I(z)$.

Define a quantity $d_f(z):=\inf \{|f(z)-w|:w \in f(\TT)\}$ which stands for the distance from $f(z)$ to the boundary of $f(\DD)$.
The next two Propositions are the Corollaries 1.4 and 1.6 of \cite{P} which are frequently used in the sequel. By the way, they are corollaries of the Koebe Distortion Theorem.

\bP\cite[Corollary 1.4]{P}\label{2503201} If $f$ maps $\DD$ conformally into $\CC$ then
\baaa{\label{2506123}\frac{1}{4}(1-|z|^2)|f'(z)|\le d_f(z)\le (1-|z|^2)|f'(z)|,}
for any $z\in\DD$.
\eP

\bP\cite[Corollary 1.6]{P}\label{2503161}
Let $a,b,c$ be positive constants and let $0<|z_0|=1-\dl<1$. If 
\baaa{\label{2504032} 0 \le 1-b\dl \le |z|\le 1-a\dl,\;\;\; |\arg z-\arg z_0|\le c\dl}
and if $f$ maps $\DD$ conformally into $\CC$ then 
\ba{|f(z)-f(z_0)|\le M_1\dl |f'(z_0)|,}
\ba{M_2^{-1}|f'(z_0)|\le |f'(z)|\le M_2|f'(z_0)|,}
where the constants $M_1$ and $M_2$ depend only on $a,b,c$.
\eP

\begin{figure}[htbp]
\includegraphics[width=0.3\textwidth]{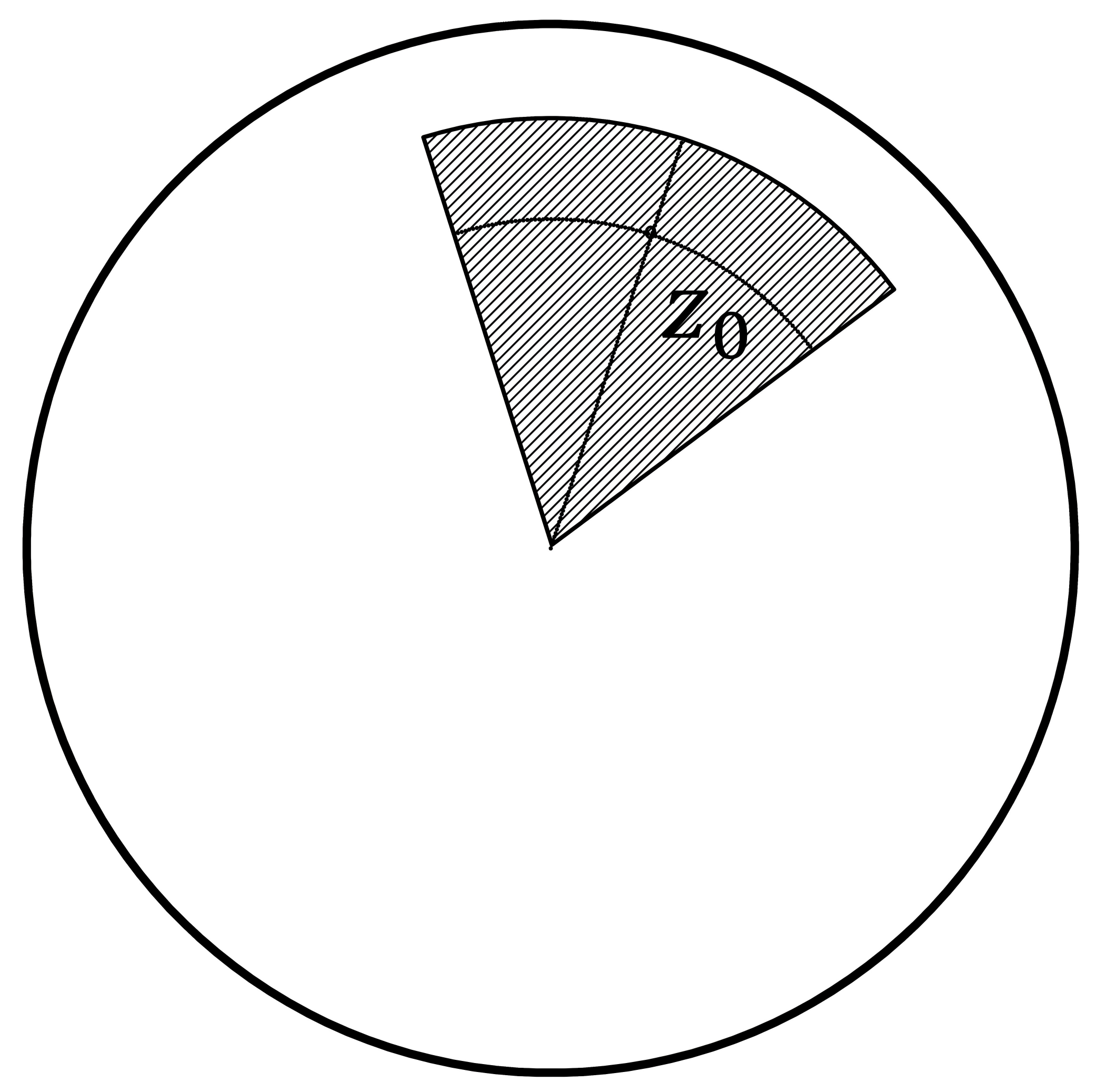}
\caption{The shaded domain is defined by (\ref{2504032})}\label{sector}
\end{figure}

Let $z_0\in \DD$. The conditions of (\ref{2504032}) describe a sectoral subdomain $D$, see Figure~\ref{sector}. For any $z\in D$, Propositions \ref{2503201} and \ref{2503161} imply 
\ba{ d_f(z_0)\ge \frac{1}{4}(1-|z_0|^2)|f'(z_0)|\ge \frac{1}{4}(1-|z_0|^2)M_2^{-1}|f'(z)| \ge \frac{1}{4M_2}\frac{1-|z_0|^2}{1-|z|^2}d_f(z).}
Interchange $z_0$ and $z$, we obtain that
\ba{ d_f(z)\ge \frac{1}{4M_2}\frac{1-|z|^2}{1-|z_0|^2}d_f(z_0).}
Then for any $z,\zt\in D$, it follows that
\ba{ d_f(z)\ge \frac{1}{16M^2_2}\frac{1-|z|^2}{1-|\zt|^2}d_f(\zt).}
Specially when $\zt=0$, it follows that 
\baaa{\label{2504033} d_f(z)\ge \frac{1}{16M^2_2}(1-|z|^2)d_f(0).}

The next proposition is also a corollary of the Koebe Distortion Theorem. However, it reveals the estimate about the ratio $|f'(z)|/|f'(w)|$ for any $z,w\in\DD$.

\bP\cite[Theorem~2.10.8]{AIM}\label{2503261}  
Suppose that $f$ is conformal in $\DD$ and $z,w\in \DD$.
Then
\ba{e^{-3\lm_{\DD}(z,w)}\le \frac{|f'(z)|}{|f'(w)|}\le e^{3\lm_{\DD}(z,w)}.}
Here $\lm_{\DD}(z,w)$ stands for the hyperbolic metric between $z$ and $w$ in $\DD$.
\eP

It is well-known the hyperbolic metric on $\DD$ is
\ba{\lm_{\DD}(z,w)=\frac{1}{2}\log\frac{1+|z-w|/|1-\bar z w|}{1-|z-w|/|1-\bar z w|},}
for any $z,w\in \DD$. Because of $\lm_{\DD}(z,w)\to \infty$ as $|w|\to 1$, we can not expect $\lm_{\DD}(z,w)$ to be bounded. However, when $|z|\ge 1/2$, let $\zt_1,\zt_2$ be the end points of $I(z)$, and let $\zt$ be the point on the non-euclidean segment from $\zt_1$ to $\zt_2$ which is nearest to $0$. Then $\lm_{\DD}(z,\zt)$ is bounded by the next corollary. A non-euclidean segment (also called a geodesic) in $\DD$ is either a circular arc that intersects the boundary $\TT$ at right angles or a diameter of $\DD$.

\bC\label{2504161}
For any $z\in \DD$ with $|z|\ge 1/2$, denote by $\zt_1,\zt_2$ the end points of $I(z)$. Let $S\sbs \DD$ be the non-euclidean segment from $\zt_1$ to $\zt_2$ and let $\zt\in S$ be the point nearest to $0$. Then $|\zt|<|z|$ and $\lm_{\DD}(z,\zt)$ is bounded.
\eC
\bpp
Without loss of generality, we assume that $z=r\ge 1/2$ then $\zt=\dl\in \RR$. By elementary geometry of planar curves $\dl=(1-\sin\pi(1-r))/\cos\pi(1-r)$. Thus it can be verified that $\dl<r$. Moreover, we obtain that
\ba{\frac{r-\dl}{1-r\dl}\le \frac{\pi-1}{\pi+1}.} 
It follows from the definition 
\ba{\lm_{\DD}(z,\zt)=\frac{1}{2}\log\frac{1+|r-\dl|/|1-r\dl|}{1-|r-\dl|/|1-r\dl|}\le \frac{1}{2}\log\pi.}
The proof is complete.
\epp

Suppose that $S$ is a non-euclidean segment in $\DD$, the next proposition provides us an estimate of the arclength of the image curve $f(S)$.

\bP\cite[Corollary 4.18]{P}\label{2503314}
If $z_0\in \DD$ and if $I$ is an arc of $\TT$ then
\ba{\ell[f(S)]<d_f(z_0)e^{[M_3/\om(z_0,I)]},}
where $S$ is the non-euclidean segment from $z_0$ to a suitable point on $I$. Moreover, $\ell[f(S)]$ is the arclength of $f(S)$ and $\om(z_0,I)$ is the harmonic measure of $I$ at $z_0$.
\eP

For a measurable subset $I$ of $\TT$, the harmonic measure of $I$ at $z$ is defined by
\ba{\om(z,I)=\frac{1}{2\pi}\int_I \frac{1-|z|^2}{|\zt-z|^2}|d\zt|.}
As $|z|\to 1$, we see that $\om(z,I)\to 0$. For a special set $I$, the next corollary provides a positive lower bounded of $\om(z,I)$.

\begin{figure}[htbp]
\includegraphics[width=0.3\textwidth]{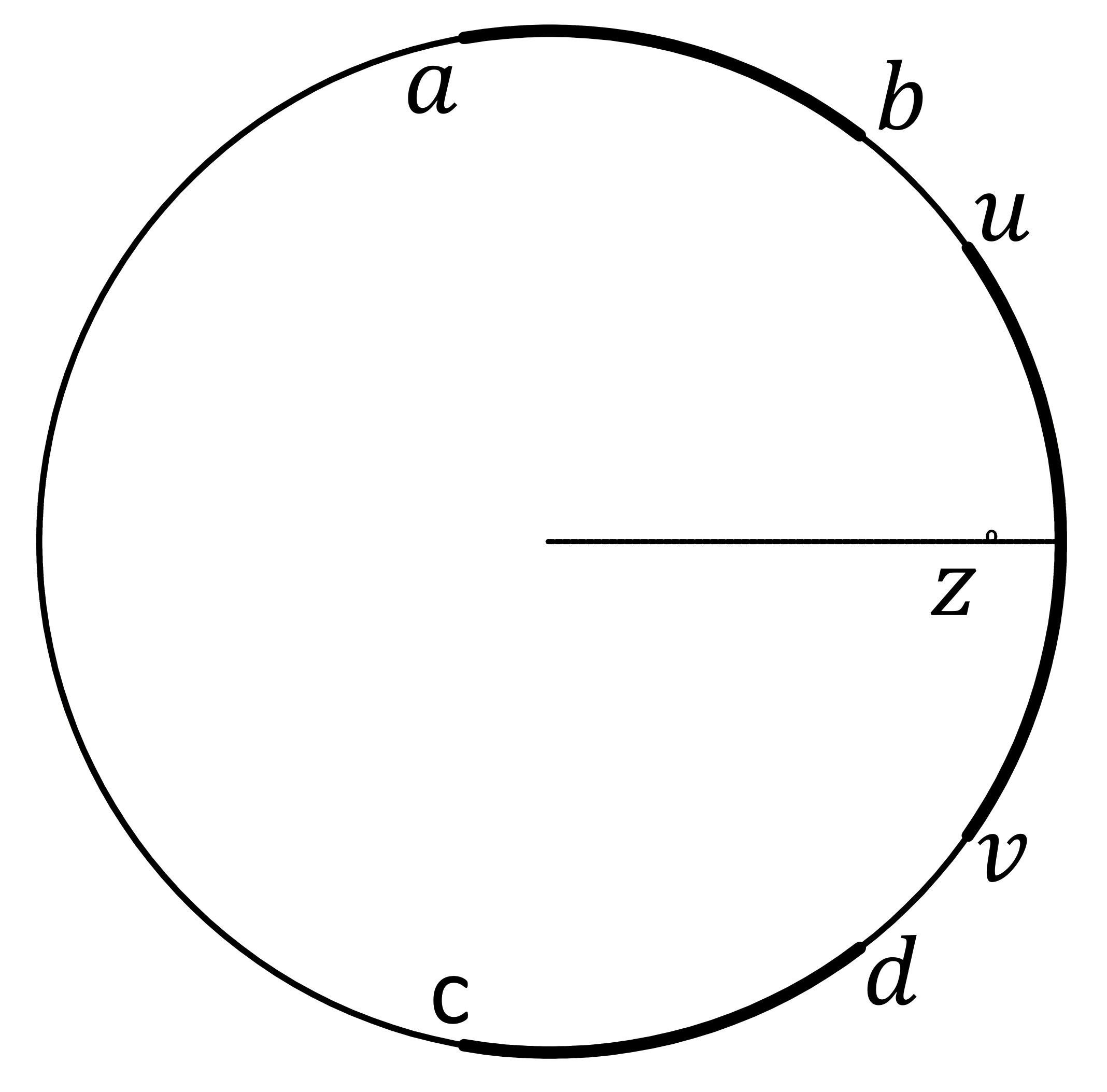}
\caption{The arcs $ab\cup cd$ is $A(z)$; the arc $uv$ is $I(z)$}\label{Az}
\end{figure}

\bC\label{2504143}
Let $re^{i\al}\in \DD$ with $r\ge 1/2$ and let $A(re^{i\al})\sbs \TT$ be the  subset defined by
\baaa{\label{2506121} A(re^{i\al}):=\{e^{i\theta}:\frac{3\pi}{2}(1-r) \le |\theta-\al|\le 2\pi(1-r)\}.}
Then the harmonic measure $\om(re^{i\al},A(re^{i\al}))$ has a positive lower bound for all $r\in [1/2,1)$.
\eC
\bpp
According to the definition of $A(re^{i\al})$, it probably consists of two subarcs which are in the complement of $I(re^{i\al})$ in $\TT$; see Figure~\ref{Az}. For any $\zt\in A(re^{i\al})$, then 
\ba{|\zt-re^{i\al}|\le 2\pi(1-r)+(1-r)<3\pi(1-r).}
Also, the measure of $A(re^{i\al})$ is 
\ba{2\times [2\pi(1-r)-\frac{3\pi}{2}(1-r)]=\pi(1-r).}  
Therefore, by the definition of harmonic measure
\ba{\om(re^{i\al},A(re^{i\al}))&=\frac{1}{2\pi}\int_{A(re^{i\al})}\frac{1-r^2}{|\zt-re^{i\al}|^2}|d\zt|\\
&\ge \frac{1}{18\pi^3}\int_{A(re^{i\al})}\frac{1-r^2}{(1-r)^2}|d\zt|\\
&=\frac{1}{18\pi^3}\frac{1+r}{1-r}\cdot \pi(1-r)\\
&\ge \frac{1}{18\pi^2}.}
The proof is complete.
\epp

The next proposition is the Gehring-Hayman theorem. It states that in simply connected planar domains,  a non-euclidean segment connecting two points is always shorter than any other curve with the same endpoints, up to a constant multiple.

\bP\cite[Theorem 4.20]{P}\label{2503264}
Let $f$ map $\DD$ conformally into $\CC$. If $z_1,z_2\in \overline{\DD}$ then
\ba{&\ell[f(S)]\le M_4\ell[f(L)],\\
& \diam f(S)\le M_5\diam f(L),}
where $S$ is the non-euclidean segment from $z_1$ to $z_2$ and $L$ is any Jordan arc in $\DD$ from $z_1$ to $z_2$.
\eP

\section{Conformal mappings onto $t$-quasidisks}\label{cmotd}

Let $\Om\sbs \CC$ be a $t$-\qdt, i.e., $\Om$ is simply connected and its boundary $\pd\Om$ is a $t$-\qct. The Riemann Mapping Theorem guarantees that there exists a conformal mapping $f:\DD\to\CC$ such that $f(\DD)=\Om$. Since $\pd\Om$ is a Jordan curve, $f$ has a continuous and injective extension to $\overline{\DD}$ by the Caratheodory Theorem (see page 18 of \cite{P}); also $f(\TT)=\pd\Om$. In this section, we study the boundary behaviour of $f$ and obtain a distortion property of $f|_\TT$ as the main result of this section. Through this section, let $N_1,\;N_2\dots$ denote suitable constants.

\bL\label{2503265}
Let $\Om$ be a $t$-\qdt. Then any two points $w_1,w_2\in \Om$ can by connected by a curve $\gm\sbs \Om$ such that 
\baaa{\label{2503262}\diam \gm\le N_1 |w_1-w_2|^t.}
\eL
\bpp
For any $w_1,w_2\in \Om$, if the line segment $[w_1,w_2]\sbs \Om$ we select $\gm=[w_1,w_2]$ and then 
\ba{\diam\gm=|w_1-w_2|\le (\diam\Om)^{1-t}|w_1-w_2|^t.}
Thus (\ref{2503262}) holds.

\begin{figure}[htbp]
\includegraphics[width=0.5\textwidth]{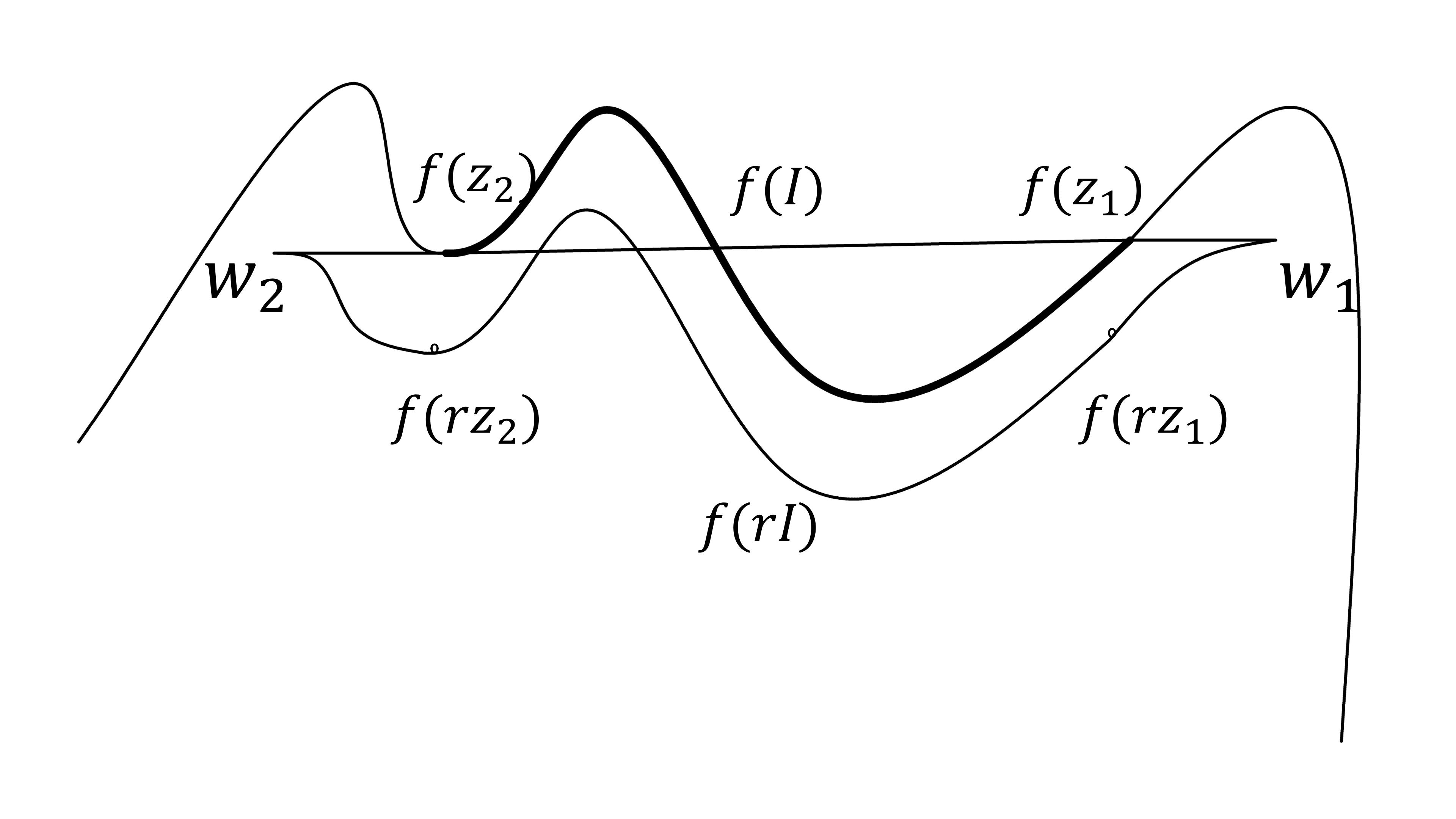}
\caption{The proof of Lemma~\ref{2503265}}\label{f5-6}
\end{figure}

Now suppose that $[w_1,w_2] \nsubseteq \Om$; see Figure~\ref{f5-6}. Let $f$ be the conformal mapping such that $f(\DD)=\Om$, and let $f(z_j)\,(\zt_j\in \TT)$ be the boundary points on $[w_1,w_2]$ nearest to $w_j$,  $j=1,2$. By the definition of $t$-\qct, one arc $I$ of $\TT$ from $z_1$ to $z_2$ satisfies
\ba{\diam f(I)\le C|f(z_1)-f(z_2)|^t.}
If $r\in (0,1)$ is close enough to $1$ then $f(rI)$ is a curve in $\Om$, and its diameter satisfies
\ba{\diam f(rI)\le 2C|f(z_1)-f(z_2)|^t<2C|w_1-w_2|^t.}
Also, we can connect $f(rI)$ with $w_j$ by curves $\gm_j\sbs \Om$ such that their arclength $\ell(\gm_j)\le |w_1-w_2|$, $j=1,2$. Let $\gm=f(rI)\cup\gm_1\cup\gm_2$. Then $\gm\sbs\Om$ is a curve from $w_1$ to $w_2$ and its diameter satisfies
\ba{ \diam \gm\le (2C+2(\diam\Om)^{1-t})|w_1-w_2|^t.}
So (\ref{2503262}) holds and the proof is complete. 
\epp

\brs
When $t=1$, Lemma~\ref{2503265} shows the linearly connected property of \qct s, which means $\Om$ does not have inward-pointing cusps (see \cite[Theorem 5.9]{TV}). However, we can not exclude the inward-pointing cusps for $t$-\qct s (see Figure~\ref{onehalf}).
\ers

\bL
Let $f$ map $\DD$ conformally onto a $t$-\qd $\Om$. For any  $z_1,z_2\in \overline{\DD}$, let $S$ be the non-euclidean segment from $z_1$ to $z_2$. Then
\baaa{\label{2503263}\diam f(S)\le N_2|f(z_1)-f(z_2)|^t.}
\eL
\bpp
By Lemma \ref{2503265}, there exists a curve $\gm$ in $\Om$ which connects $f(z_1)$ and $f(z_2)$ such that 
\ba{\diam \gm\le N_1|f(z_1)-f(z_2)|^t.} 
Let $L:=f^{-1}(\gm)$. Then $L\sbs\DD$ and connects $z_1,z_2$. It follows from Proposition \ref{2503264} that
\ba{\diam f(S)\le M_5\diam f(L) \le  N_1M_5|f(z_1)-f(z_2)|^t.} 
The proof is complete.
\epp

Recall that for any $re^{i\al}\in\DD$, we have defined the subsets
\ba{B(re^{i\al})=\{\rho e^{i\theta}: r\le \rho\le 1,\; |\theta-\al|\le \pi(1-r) \},\;\; I(re^{i\al})=B(re^{it})\cap \TT.}
Also, $d_f(re^{i\al})=\inf \{|f(re^{i\al})-w|:w \in f(\TT)\}$ which is the distance from $f(re^{i\al})$ to the boundary of $f(\DD)$. The next lemma shows that the diameter of $f(B(z))$ is dominated by $(d_f(z))^t$.
 
\bL\label{2503171} Let $f$ map $\DD$ conformally onto a $t$-\qd $\Om$, and let $c>0$ be fixed such that $d_f(0)\ge c\diam\Om$. Then 
\baaa{\label{2506122}\diam f(B(z))\le N_3(d_f(z))^t, \;\;\; \forall z\in \DD.}
\eL
\bpp
Let $\Gm:=\pd\Om$. Then $\Gm$ is a $t$-\qct. For any $z=re^{i\al}\in \DD$, if $0\le r\le 1/2$, we know that  
\baaa{\label{2503312}\diam f(B(z))\le \diam \Om = (\diam \Om)^{1-t}(\diam \Om)^t\le (\diam \Om)^{1-t} c^{-t}(d_f(0))^t.}
Since $z$ and $0$ are contained in the sectoral subdomain $D$ which is defined by $1/2$ in the form of (\ref{2504032}). It follows from (\ref{2504033}) that
\baaa{\label{2503313} d_f(0)\le  \frac{16M^2_2}{1-r^2}  d_f(z) \le \frac{64M^2_2}{3}d_f(z).}
Combining (\ref{2503312}) and (\ref{2503313}) yields that 
\baaa{\label{2503061} \diam f(B(z)) \le (\diam \Om)^{1-t}c^{-t}(\frac{64M^2_2}{3})^t(d_f(z))^t.}

If $1/2<r<1$. Let $\zt_1$ and $\zt_2$ be the end points of $I(z)$. Let $I^\pm\sbs \TT$ be the subsets defined by (\ref{2506121}), i.e., they are subarcs consisting of all points whose argument $\theta$ satisfies
\ba{\frac{3\pi}{2}(1-r)\le |\theta-\al|\le 2\pi(1-r).}
Then $I^\pm\cap I(z)=\emptyset$, meanwhile $I^\pm$ do not intersect the segment $S_0:=[-e^{i\al},0]$. Suppose that $I^+$ is the subarc between $\zt_1$ and $-e^{i\al}$, while $I^-$ is the subarc between $\zt_2$ and $-e^{i\al}$. By Proposition \ref{2503314} and Corollary \ref{2504143}, there are non-euclidean segments $S^\pm$ from $z$ to $I^\pm$ such that 
\baaa{\label{2503064}\ell[f(S^+\cup S^-)]\le N_4 d_f(z).}
Here $N_4=2\exp(18\pi^2M_3)$. Let $L:=f(S^+\cup S^-)$, and let the two end points of $L$ be $a$ and $b$. Let $H$ and $G$ be the two  components of $\Om\setminus L$, and let $\Gm(a,b)\cup L$ be the boundary of $H$. It follows that 
\baaa{\label{2504081} \diam H &\le \diam L+\diam \Gm(a,b)\le \diam L+C|a-b|^t\\
& \le \diam L+C(\diam L)^t\le ((\diam \Om)^{1-t}+C)(\diam L)^t. }
Here $C$ is the constant in the definition of $t$-\qct. 

\begin{figure}[htbp]
\includegraphics[width=0.3\textwidth]{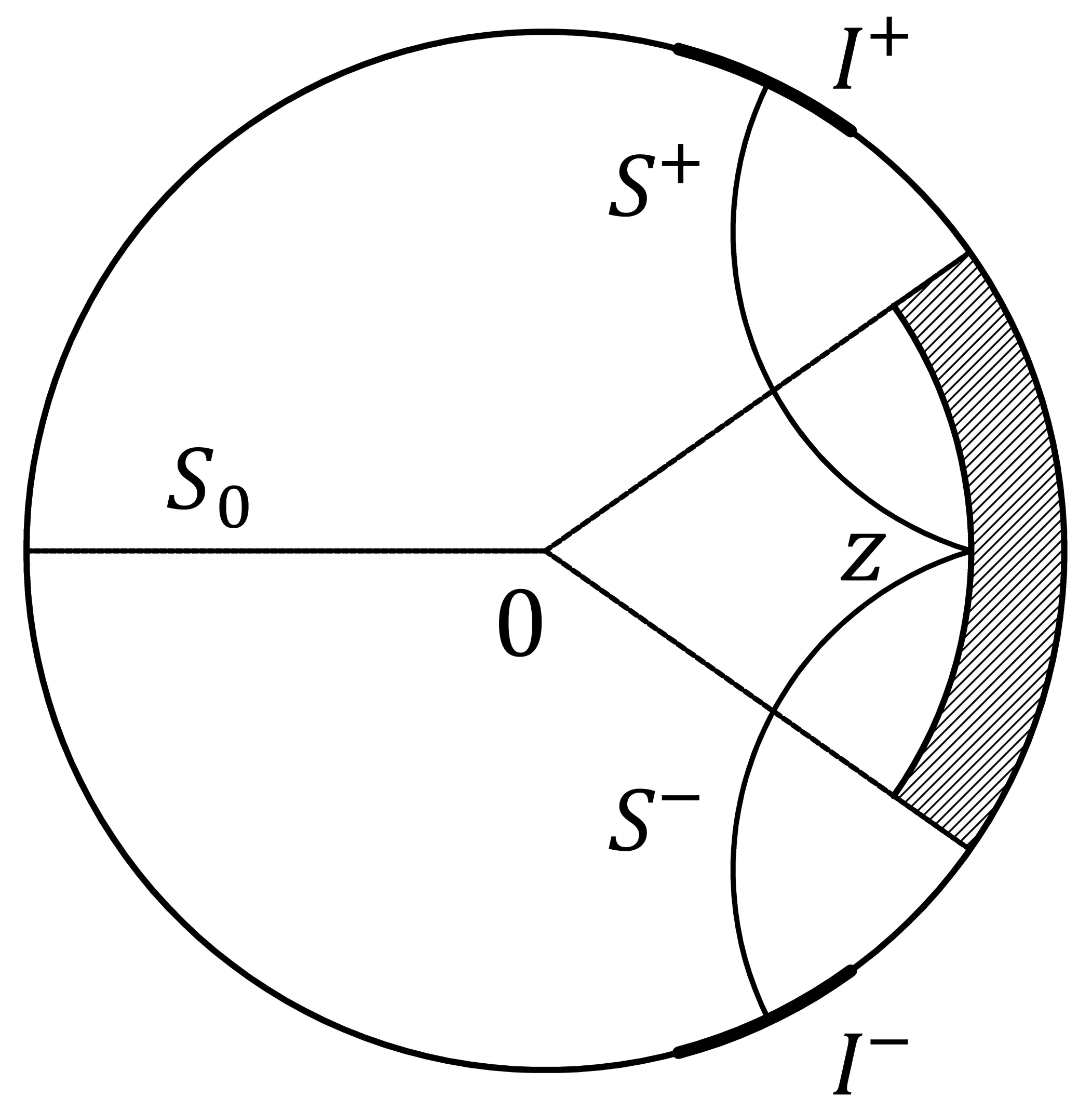}
\caption{The shaded domain is $B(z)$}\label{Spm}
\end{figure}

Our choice of $I^\pm$ implies that $S^+\cup S^-$ separates $B(z)$ from $S_0$;  see Figure~\ref{Spm}, then $L$ separates $f(B(z))$ and $f(S_0)$. If $f(B(z))\sbs \overline{H}$, it follows from (\ref{2504081}) that 
 \baaa{\label{2503063}\diam f(B(z))\le\diam H\le N_5 (\diam L)^t.}
By (\ref{2503064}) and (\ref{2503063}), we obtain that
\ba{\diam f(B(z))\le N_5(\diam L)^t\le N_5(\ell(L))^t\le N_5N_4^t(d_f(z))^t.}
If $f(B(z))\sbs \overline{G}$ then $f(S_0)\sbs H$, and then
\ba{d_f(0)\le \diam f(S_0)\le \diam H\le  N_5N_4^t(d_f(z))^t.}
Since
\ba{\diam f(B(z))\le \diam G\le \diam \Om \le cd_f(0),}
it follows that 
\ba{\diam f(B(z))\le c N_5N_4^t(d_f(z))^t.}
Thus we conclude that 
\baaa{\label{2503062}\diam f(B(z))\le (c+1) N_5N_4^t(d_f(z))^t.}
whenever $1/2< r <1$. 

Combining (\ref{2503061}) and (\ref{2503062}) implies (\ref{2506122}). The proof is complete. 
\epp

For any $z\in\DD$, let $\zt_1,\zt_2$ be the endpoints of $I(z)$. Then
\ba{|f(\zt_1)-f(\zt_2)|\le \diam f(B(z)),}
and then (\ref{2506123}) and (\ref{2506122}) imply that
\ba{|f(\zt_1)-f(\zt_2)|\le  N_3(d_f(z))^t\le 2^tN_3(1-|z|)^t|f'(z)|^t.}
We have derived the following lemma.

\bL\label{2504104}
Let $f$ map $\DD$ conformally onto a $t$-\qdt, and let $z\in \DD$. Suppose that $\zt_1,\zt_2$ are the endpoints of $I(z)$, then there exists a constant $M>0$ such that 
\ba{|f(\zt_1)-f(\zt_2)|\le  M(1-|z|)^t|f'(z)|^t.}
\eL

Lemma~\ref{2504104} established the upper bound of $|f(\zt_1)-f(\zt_2)|$,  the next lemma provides us the lower bound.

\bL\label{2504103}
Let $f$ map $\DD$ conformally onto a $t$-\qdt, and let $z\in \DD$ with $|z|\ge 1/2$. Then
\ba{|f(\zt_1)-f(\zt_2)|^t\ge N (1-|z|)|f'(z)|,}
where $\zt_1, \zt_2$ are the endpoints of $I(z)$ and $N$ is a suitable constant.
\eL
\bpp
Let $r=|z|$, $\zt_j=e^{i\theta_j}$ $j=1,2$, and let $S$ be the non-euclidean segment from $\zt_1$ to $\zt_2$ in $\DD$. Select $\zt\in \TT$ such that $\dl \zt$ be the point on $S$ that is nearest to $0$, then $r\zt=z$. Since $r \ge 1/2$, Corollary \ref{2504161} implies that $\dl\le r$ and that
\ba{\lm_{\DD}(r\zt,\dl\zt)\le \frac{1}{2}\log\pi.}
By Proposition \ref{2503261}, we obtain that
\ba{\frac{1}{\pi^{3/2}}|f'(\dl\zt)|\le |f'(r\zt)|\le \pi^{3/2}|f'(\dl\zt)|.}
 It follows from (\ref{2506123}) that
\baaa{\label{2504102} d_f(\dl \zt)\ge \frac{1}{4}(1-\dl)|f'(\dl\zt)|\ge \frac{1}{4\pi^{3/2}}(1-\dl)|f'(r\zt)|\ge \frac{1}{4\pi^{3/2}}\frac{1-\dl}{1-r}  d_f(r\zt).}
Then we conclude that
\ba{N_2|f(\zt_1)-f(\zt_2)|^t\ge \diam f(S)\ge d_f(\dl\zt)\ge \frac{1}{4\pi^{3/2}}d_f(z),}
by (\ref{2503263}) and (\ref{2504102}).
Therefore,
\ba{|f(\zt_1)-f(\zt_2)|^t\ge \frac{1}{4\pi^{3/2}N_2}d_f(z)\ge  \frac{1}{4\pi^{3/2}N_2}\frac{1}{4}(1-|z|)|f'(z)|.}
The last inequality follows from (\ref{2506123}).
The proof is complete.
\epp

In the beginning of this section, we mentioned that if $f$ maps $\DD$ conformally onto a $t$-\qd $\Om$, it could be extended onto $\DD\cup\TT$. We use the same notation $f$ to denote the induced boundary value $f:\TT\to \pd\Om$. The next theorem shows that $f|_\TT$ is ``almost" $(\rho,t^2)$-\qst.

\bT\label{2505274}
Let $\Om$ be a $t$-\qd and let $f$ be the conformal mapping with the induced boundary value $f:\TT\to \pd\Om$. Then there exists an increasing injection $\eta:[0,+\infty)\to [0,+\infty)$ such that
\baaa{\label{2505262}\frac{|f(a)-f(x)|}{|f(b)-f(x)|^{t^2}}\le \eta(k)}
for any $a,b,x\in \TT$ satisfy $|a-x|=k|x-b|$ with $k\ge 1$. 
\eT
\bpp
Without loss of generality, we assume that $x=1$. Let $a=e^{i\theta_1}$ and $b=e^{i\theta_2}$ where $\theta_1,\theta_2\in (-\pi,\pi]$. Since $|a-1|\ge|b-1|$, there are $\zt_1=r_1e^{i\al}$ and $\zt_2=r_2e^{i\bt}$ with $1/2<r_1\le r_2<1$ and $|\al|=\pi(1-r_1)$, $|\bt|=\pi(1-r_2)$ such that
\ba{I(\zt_1)=\TT(a,1)\;\;\text{ and }\;\; I(\zt_2)=\TT(b,1).}
It is not hard to see that the arclength of $I(\zt_1)$ is $2\pi(1-r_1)$ while the arclength of $I(\zt_2)$ is $2\pi(1-r_2)$. It is well-known that the unit circle is a chord-arc curve with ratio constant $\pi/2$, i.e., for any two points on $\TT$, the length of the shortest arc between them is comparable to the length of the chord connecting them. Then $2\pi(1-r_1)\le \pi/2|a-1|$. It follows that
\ba{4(1-r_1)\le |a-1|\le 2\pi (1-r_1).}
Similarly,
\ba{4(1-r_2)\le |b-1|\le 2\pi (1-r_2).}
Then we obtain that
\baaa{\label{2506042} \frac{2}{\pi}\frac{1-r_1}{1-r_2}\le \frac{|a-1|}{|b-1|}\le \frac{\pi}{2}\frac{1-r_1}{1-r_2}.}
Lemma \ref{2504104} implies that
\ba{|f(a)-f(1)|\le M(1-|\zt_1|)^t|f'(\zt_1)|^t,}
while Lemma \ref{2504103} implies that
\ba{|f(1)-f(b)|^t\ge N (1-|\zt_2|)|f'(\zt_2)|.}
It follows from (\ref{2506042}) that
\baaa{\label{2505261}\frac{|f(a)-f(1)|}{|f(1)-f(b)|^{t^2}}\le \frac{M}{N^t}\frac{(1-r_1)^t}{(1-r_2)^t}\frac{|f'(\zt_1)|^t}{|f'(\zt_2)|^t}\le \frac{M}{N^t}(\frac{k\pi}{2})^t \frac{|f'(\zt_1)|^t}{|f'(\zt_2)|^t}.}
Since $1-r_1\ge 1-r_2$, then $|\al-\bt|\le \pi(1-r_1)+\pi(1-r_2)\le 2\pi(1-r_1)$. It follows from (\ref{2505241}) and (\ref{2506042}) that
\ba{\log\frac{|f'(\zt_1)|}{|f'(\zt_2)|}\le & \frac{4r_2}{1-r_2^2}|\al-\bt|+\log\frac{1-r_1}{1-r_2}+\log 8\\
\le & \frac{8\pi r_2(1-r_1)}{1-r_2^2}+\log(4\pi \frac{|a-1|}{|b-1|})  \\
\le & 4\pi^2 \frac{|a-1|}{|b-1|}+\log(4\pi \frac{|a-1|}{|b-1|})\\
= & 4\pi^2 k+\log(4\pi k),}
thus
\ba{\frac{|f'(\zt_1)|}{|f'(\zt_2)|}\le 4\pi k e^{4\pi^2k}.}
Therefore we obtain from (\ref{2505261}) that 
\ba{\frac{|f(a)-f(1)|}{|f(b)-f(1)|^{t^2}}\le (4\pi k e^{4\pi^2k})^t\frac{M}{N^t}(\frac{k\pi}{2})^t=: \eta (k).}
It is not hard to verify that the function $\eta$ is an increasing injection from $[0,\infty)$ onto itself. The proof is complete.
\epp

So far, we can not say that the boundary value $f|_\TT$ is $(\eta,t^2)$-\qst, since (\ref{2505262}) does not hold for $k<1$. However, when $|a-x|=|b-x|$ i.e., $k=1$, we obtain that
\ba{\frac{|f(a)-f(x)|^{1/t}}{|f(x)-f(b)|^t}\le \eta(1) .}
This could be regarded as a generalization of Ahfors' $M$-condition or $M$-\qst, which is satisfied by the boundary values of quasiconformal mappings on $\DD$ (see \cite{A}). It implies that the boundary value $f|_\TT$ is weakly $(R,t^3)$-\qs from the following corollary.

\bC\label{2512033}
Let $\Om$ be a $t$-\qd and let $f$ be the conformal mapping with the induced boundary value $f:\TT\to \pd\Om$. Then $f$ is weakly $(R,t^3)$-\qs on $\TT$.
\eC
\bpp
Let $a,b,x\in \TT$ be arbitrary three distinct points such that $|a-x|\le |b-x|$. Then $f(a),f(b),f(x)\in \Gm$, where $\Gm=\pd\Om$ which is a $t$-\qct, and then
\ba{\diam \Gm(f(x),f(b))\le C|f(b)-f(x)|^t}
for some positive constant $C$. Here $\Gm(f(x),f(b))$ is the subarc of $\Gm$ between $f(x)$ and $f(b)$ which has a smaller diameter. 

Suppose that $\Gm(f(x),f(b))=f(\TT(x,b))$. Here $\TT(x,b)$ is the subarc of $\TT$ between $x$ and $b$ which has a smaller diameter. If $f(a)\in \Gm(f(x),f(b))$ i.e., $a\in \TT(x,b)$, then $|f(a)-f(x)|\le \diam \Gm(f(x),f(b))$ and then 
\baaa{\label{2512031}|f(a)-f(x)|\le C|f(b)-f(x)|^t.} 
If $a\notin \TT(x,b)$, we can select a point $c\in \TT(x,b)$ such that $|a-x|=|c-x|$ then we obtain that
\ba{|f(c)-f(x)|\le C|f(b)-f(x)|^t.}
By Theorem \ref{2505274}, it follows that 
\ba{\frac{|f(a)-f(x)|}{|f(c)-f(x)|^{t^2}}\le \eta(1).}
So we conclude that
\baaa{\label{2512032}|f(a)-f(x)|\le  C^{t^2}\eta(1)|f(b)-f(x)|^{t^3}.}

If $\Gm(f(x),f(b))=f(\TT\setminus\TT(x,b))$, the proof is similar to the case of $\Gm(f(x),f(b))=f(\TT(x,b))$, so we omit the details. 

Therefore (\ref{2512031}) and (\ref{2512032}) implies that $f$ is weakly $(R, t^3)$-\qst. The proof is complete.
\epp

\brs
According to the definition of weakly $(R,t)$-\qs mappings, the constant $R$ is required to be positive. However, we may assume that $R>1$ since
\ba{|f(a)-f(x)|\le R|f(b)-f(x)|^t\le \frac{2}{R}\cdot R|f(b)-f(x)|^t}
whenever $R<1$. 
\ers

We conclude this section with the following theorem, which is the main result. By letting $t=1$, we can replicate the Quasicircle Theorem from our result.

\bT\label{2512034}
Let $\Om$ be a $t$-\qd and let $f$ be the conformal mapping with the induced boundary value $f:\TT\to \pd\Om$. Then $f$ is $(\rho,t^6)$-\qs on $\TT$.
\eT
\bpp
Without loss of generality, we may assume that $\diam\Om=1$, as its $t$-\qd property is invariant under scaling. For any $k\in(0,\infty)$, let $a,b,x\in \TT$ be distinct points with $|a-x|=k|b-x|$. Set $m:=|f(a)-f(x)|/|f(b)-f(x)|^{t^6}$. By Theorem \ref{2505274}, there exists an increasing injection $\eta:[0,\infty)\to[0,\infty)$ such that
$|f(a)-f(x)|/|f(b)-f(x)|^{t^2}\le \eta(k)$ when $k\ge 1$. It follows that $m\le \eta(k)|f(b)-f(x)|^{t^2-t^6}\le (\diam\Om)^{t^2-t^6}\eta(k)=\eta(k)$ when $k\ge 1$.  

The Corollary~\ref{2512033} states that $f$ is weakly $(R,t^3)$-\qs on $\TT$. Now if $k<1$ i.e., $|a-x|<|b-x|$, and then the weakly $(R,t^3)$-\qs property of $f$ implies that $|f(a)-f(x)|\le R|f(b)-f(x)|^{t^3}$. So $m\le R(\diam\Om)^{t^3-t^6}=R$ when $k<1$. We select by induction a sequence of points $b_0,b_1,\dots$ of $\TT$ such that $b_0=b$ and that $|b_{j-1}-x|/4<|b_j-x|<|b_{j-1}-x|/3$ for all $j\ge 1$. Then $b_j\to x$ as $j\to\infty$. Let $s$ be the least integer $j$ such that $|b_j-x|<|a-x|$. Since $|b_j-x|>|b-x|/4^{j}$, then $|a-x|>|b-x|/4^s$, i.e., $4^{-s}<k$.

Suppose that $0\le i<j<s$. Then $|x-b_j|\le |x-b_{i+1}|<|x-b_i|/3\le |x-b_j|/3+|b_i-b_j|/3$, which implies that $2|x-b_j|\le |b_i-b_j|$. Because of $|a-x|\le |b_j-x|$, it implies that 
\ba{|a-b_j|\le |a-x|+|x-b_j| \le |b_i-b_j|.} 
Thus 
\ba{|f(a)-f(x)|\le |f(a)-f(b_j)|+|f(b_j)-f(x)|\le 2R|f(b_i)-f(b_j)|^{t^3}.}
Thus the distance between the points $f(b_j)$, $0\le j\le s-1$, are at least $|f(a)-f(x)|^{1/{t^3}}/(2R)^{1/{t^3}}$. On the other hand, $|b_j-x|\le |b-x|$ implies that $|f(b_j)-f(x)|\le R|f(b)-f(x)|^{t^3}$. Since any subset of $\RR^2$ is homogeneously totally bounded with $\vp(\al)=4(\al\sqrt{2}+1)^2$ (see the second paragraph after Proposition~\ref{2503141}), and since $m\le R$ when $k<1$, it follows that $m<2R(2R)^{t^3}$. Thus 
\ba{\al:=\frac{R(2R)^{1/t^3}}{m^{1/t^3}} >\frac{1}{2}.} 
Therefore, the homogeneously totally bounded property implies that
\ba{s\le \vp( \frac{R|f(b)-f(x)|^{t^3}}{|f(a)-f(x)|^{1/{t^3}}/(2R)^{1/{t^3}}})=\vp(\al).}
Since $4^{-s}<k$, we obtain that 
\ba{\log\frac{1}{k}<\vp(\frac{(2R)^{1/{t^3}}R}{m^{1/{t^3}}})\log 4.}
Directly calculation shows that when $k< 1/256$
\ba{m^{1/{t^3}}\le \{[\frac{(\log(1/k))}{\log 256}]^{1/2}-1\}^{-1}\sqrt{2}R(2R)^{1/{t^3}}=:\psi(k)^{1/t^3}.}
It is obvious that $\psi(k)\to 0$ as $k\to 0$. Then $\psi$ can be extended to  an increasing continuous function on $[0,1/512]$ such that $\psi(0)=0$. 

In $[1/512,1]$, we can define an increasing linear function $p(k)$ such that $p(1/512)=K_1\psi(1/512)$ and $p(1)=K_2\eta(1)$. Here $K_1$ and $K_2$ are suitable numbers such that $m\le p(k)$ for $k\in [1/512,1]$ and that $p(k)$ is increasing. These $K_1$ and $K_2$ exist since $m\le R$, which means that $m$ is bounded when $k\in (0,1)$. So we can define an increasing injection $\rho:[0,\infty)\to [0,\infty)$ by 
\ba{ \rho(k)=\begin{cases}
K_1\psi(k), & 0\le k\le 1/512 ;\\
p(k), & 1/512<k<1;\\
K_2\eta(k), & k\ge 1.
\end{cases}
}
Therefore $f$ is $(\rho,t^6)$-\qs on $\TT$, and the proof is complete.
\epp

\brs
We consider that, in this theorem, the exponent $6$ of $t^6$ may not be optimal. Our goal is to achieve the exponent $1$, i.e., the boundary value $f$ is $(\rho,t)$-\qs on $\TT$. However, to attain this objective, it might require adopting some other approaches, which we will continue to explore in future work. While this result is less than ideal, we have successfully generalized the Quasicircle Theorem. 
\ers

\section*{Acknowledgments }

The author thanks the anonymous reviewers for their helpful and
valuable suggestions on improving this paper.

\section*{Data Availability }
Data sharing is not applicable to this article as no datasets were generated or analysed during the current study.

\section*{Conflict of interest}

The author declares no conflicts of interest.

\end{document}